# Multi-objective optimization of double curvature arch dams subjected to seismic loading using charged system search

S. Talatahari[1,2]* M.T. Aalami[1] and R. Parsiavash[1]

[1] *Department of Civil Engineering, University of Tabriz, Tabriz, Iran*

[2] *Faculty of Engineering & IT, University of Technology Sydney, Ultimo, NSW 2007, Australia*

**Abstract**

This paper presents a multi-objective formulation for optimization of arch dams. The objective is to simultaneously minimize the concrete volume and stress state of the dam body. Instead of a single design, the multi-objective problem provides a set of Pareto solutions in which a decision maker approach can select best answer according to its preferences. To solve this new problem, a multi-objective charged system search as an efficient meta-heuristic is also developed. The proposed algorithm is evaluated on a real-size double curvature arch dam. The computational results show that the proposed methodology is highly competitive with existing arch dam optimization algorithms.

*Corresponding author: *E-mail address*: Siamak.talat@gmail.com (S. Talatahari)

# 1. Introduction

Economy and safety of arch dams especially for their seismic performance are of high importance. Earthquake acceleration produces not only dynamic loads on the dam body but also hydrodynamic pressures in reservoir, which intensifies resultant loads to the arch dam due to fluid-structure interaction effects. In order to reduce the complexity of problem in some early researches, seismic loads were ignored [1-4] and in some others, hydrodynamic pressures and dam-reservoir interaction effects were simplified [5, 6]. In recent years due to technological advancements in structural modeling and processing approaches, there is an increasing interest in optimal design of arch dams considering the dynamic analysis and fluid-structure interactions. Employing appropriate structural analyzing approaches as well as advanced optimization techniques are unavoidable to fulfill this purpose.

Nowadays, finding an efficient and reliable methodology using new optimization techniques has become so necessary for the structural design problems. In this field, some successful applications of optimization algorithms have been reported in the literature [7-12]. However, in the field of dam optimization, many of them used mathematical programming approaches. For effective implementation of these methods, the variables and cost function need to be continuous. Furthermore, a good starting point is vital for these methods to be executed successfully. In many optimization problems, prohibited zones, side limits, and non-smooth or non-convex cost functions need to be considered. As a result, these non-convex optimization problems cannot be solved by the traditional mathematical programming methods. Therefore, meta-heuristic approaches as the alternative methods were introduced and applied to overcome drawbacks of conventional approaches. Recent developments in this field have led to a renewed interest in employing them on the structural optimization field. The Charged system search method (CSS) proposed by Kaveh and Talatahari [13] is one of these methods that attracted many attentions especially in the structural optimization [14-18].

On the other hand, the arch dam optimization problem can be formulated as a multiple criteria decision-making problem in which optimal decisions need to be taken in the presence of trade-offs between two or more conflicting objectives. Therefore, there exist a set of optimal solutions instead of only a single solution. Until now, all the researches formulated the arch dam optimization as a one-objective problem. In most of these researches, the concrete volume is considered as the main objective function and some behavior, geometry or stability conditions as constraints [12, 19]. In this study, an efficient methodology is proposed for shape optimal design of double-curvature arch dams subjected to earthquake loading considering multi-objective optimization concept. The CSS method is then modified and employed to meet this new defined problem. The time history analysis of arch dam models is performed under the gravity load, hydrostatic/hydrodynamic pressures and earthquake loads. The concrete volume and stress state of dam body are considered as the objective functions and shape parameters of the dam are used as design variables. Some geometry and stability constraints are also defined to ensure construction requirements meet.

The rest of the paper continuous by the geometric and finite element modeling of the double curvature arch dam in Sections 2 and 3, respectively. It will then go on to the problem statements in Section 4. In the next two sections, basic aspects and the characteristics of the standard and multi-objective CSS are explained. Then, the multi criteria decision making is presented in Section 7. The case study is investigated in Section 8 to verify the efficiency of the proposed method and finally some concluding remarks are provided in the last section.

## 2. Dam Geometry

The double-curvature arch dams feature horizontal as well as vertical arches. Therefore, both curvature and thickness of the dam body change in both directions. To define the geometrical model of this structure, shape of crown cantilever should be determined at first. In this study, one polynomial of second order is used to determine the curve of upstream boundary as shown in Fig. 1. [20, 21].

$$g(z) = \frac{\gamma_s z^2}{2\beta h} - \gamma_s z \qquad (1)$$

where $\gamma_s = \cot\alpha$ is slope of overhang at crest, h is the dam height and the point where the slope of the upstream face equals to zero is $z = \beta h$. By dividing the dam height into n segments (n+1 controlling levels) and specifying the thickness of each level, a polynomial is employed to determine the thickness of the crown cantilever.

$$t_c(z) = \sum_{i=1}^{n+1} L_i(z) t_{ci} \tag{2}$$

where $L_i(z)$ is a Lagrange interpolation formula and $t_{ci}$ is the thickness of central vertical section at $i$th level. So, the thickness of dam in other points can be obtained by interpolating the values of thickness at controlling levels. $L_i(z)$ can be expressed as:

$$L_i(z) = \prod_{\substack{m=1 \\ m \neq i}}^{n+1} \frac{z - z_m}{z_i - z_m} \tag{3}$$

in which, $z_i$ denotes the z coordinate of $i$th level.

As shown in Fig. 2, for the purpose of symmetrical canyon and arch thickening from crown to abutment, the shape of the horizontal section of a parabolic arch dam is determined by the following two parabolas [20, 21].

$$y_u(x,z) = \frac{x^2}{2r_u(z)} + g(z) \tag{4}$$

$$y_d(x,z) = \frac{x^2}{2r_d(z)} + g(z) + t_c(z) \tag{5}$$

where $y_u(x,z)$ and $y_d(x,z)$ are parabolas of the upstream and downstream faces, respectively. $r_u$ and $r_d$ are the radius of curvature of upstream and downstream curves at z direction and can be interpolated by $L_i(z)$ according to following equations of $n$th order:

$$r_u(z) = \sum_{i=1}^{n+1} L_i(z) r_{ui} \tag{6}$$

$$r_d(z) = \sum_{i=1}^{n+1} L_i(z) r_{di} \tag{7}$$

in which $r_{ui}$ and $r_{di}$ are the values of $r_u$ and $r_d$ at the controlling levels, respectively.

## 3. Finite Element Modeling of Arch Dams

A finite element model utilizing two hundred and eighty 8-node solid elements including 495 nodes for the dam body and one thousand and four hundred 8-node fluid element including 2145 nodes for the reservoir has been developed as shown in Fig. 3. This number of nodes and elements are not constant during the optimization process due to changes in dimensions of the arch dam and mesh generation in every analysis, so it will be varied when needed. Solid elements have three degrees of freedom at each node including translations in the nodal x, y and z directions. Fluid elements have four degree of freedom per node including translations in the x, y and z direction and pressure. The translation, however, are applicable only at nodes that are in the interface. To prevent extra complexities, only the dam and reservoir interaction are considered and foundation rock is assumed rigid. The dam is treated as a 3D linear structure. Fluid medium is assumed as homogeneous, isotropic, irrotational and inviscid by linear compressibility. Fluid-structure interaction (FSI) effects are taken to account between reservoir and dam body as well as reservoir walls (Fig. 4). Effects of the surface wave are ignored. In order to consider the damping effect arising from the propagation of pressure waves in the upstream direction, instead of a Sommerfield-type radiation boundary condition, the reservoir length is selected as three times the reservoir depth, and zero pressure is imposed on all nodes of the far end boundary as recommended by [22, 23]. In this research, three loading cases including gravity, static and seismic loading have been considered. The dynamic analysis of the system is performed using the Newmark time integration method.

### 3.1. Verification of the analytical model

In order to verify analytical model, a finite element model (FEM) of the morrow point arch dam-reservoir system has been developed (Fig. 3) This double curvature thin-arch concrete structure is located in the Gunnison River in west-central Colorado about 35 km east of Montrose. The dam has 142.65 m high and 220.68 m long along the crest and its thickness varies from 3.66 m at the crest to 15.85 m at the base level [24]. The dam body construction required 273600 $m^3$ of concrete. Material properties of both the dam body and water are presented in Table1. The natural frequencies of arch dam with both empty and full reservoir are obtained using finite element analysis (FEA) and compared to experimental and analytical results from other works [25, 26] summarized in Table2. It can be observed that good correspondence has been achieved between the results of present work with those reported in the literatures.

## 4. Formulation of Arch Dam Multi-Objective Optimization

Multi-objective optimization is a framework in which many objective functions (often two) are desired to have equal treatments. The solution of this problem is a set of multiple sub-solutions, which optimizes simultaneously the objectives. In addition, the feasible solutions have to comply with certain constraint conditions. Eventually, the decision is acquired after the optimization process has finished. This problem can be formulated as follows:

*Minimize:* $\quad (fit_1(\mathbf{X}), fit_2(\mathbf{X}), ..., fit_k(\mathbf{X}), ..., fit_N(\mathbf{X}))$ (8)

*Subject to:* $\quad \varnothing_i(X) > 0 \quad (i = 1, 2, 3, ..., p)$ (9)

where $fit_k(X)$ and $\varnothing_i(X)$ denote the multiple different objectives and the constraint condition respectively. *N and p* are the number of objectives and constraints, respectively. **X** is the feasible set of decision variables. In this paper, two objectives are considered.

Objective 1: Concrete volume of the arch dam. It can be calculated by integrating the dam surfaces as:

*Minimize:* $\quad fit_1(X) = \iint_A |y_u(x,z) - y_d(x,z)| dA$ (10)

Objective 2: Failure criterion function of Willam and Warnke. For concrete structures, it is defined as [27]:

*Minimize:* $\quad fit_2(X) = \left(\dfrac{F}{f_c} - \dfrac{S}{S_f}\right)_{n,t} \quad n = 1, 2, ..., n_d \quad t = 1, 2, ..., T$ (11)

where, *F* is the function of principal stress state ($\sigma_{xp}, \sigma_{yp}, \sigma_{zp}$) in which, $\sigma_{xp}, \sigma_{yp}$ *and* $\sigma_{zp}$ are principal stresses in principal directions *x, y and z* respectively. T is the earthquake duration. $n_d$ is the total number of nodes in the finite element model. $s_f$ is safety factor which for the earthquake loading may be chosen as $s_f = 1$ [28]. S is failure surface expressed in terms of principal stresses. A total of five input strength parameters are needed to define the failure surface as well as an ambient hydrostatic stress state: ultimate uniaxial compressive strength of concrete $(f_c)$, ultimate uniaxial tensile strength of concrete $(f_t)$, ultimate biaxial compressive strength of concrete $(f_{cb})$, ambient hydrostatic stress state $\sigma_h^a$, ultimate compressive strength for a state of biaxial compression superimposed on hydrostatic stress state $(f_1)$ and ultimate compressive strength for a state of uniaxial compression superimposed on hydrostatic stress state $(f_2)$. However, the failure surface can be specified with a minimum of two constants, $f_t$ and $f_c$. The other three constants can be obtained using Eqs. 12-14 according to willam and Warnke [27]:

$$f_{cb} = 1.2 f_c \tag{12}$$

$$f_1 = 1.45 f_c \tag{13}$$

$$f_2 = 1.725 f_c \tag{14}$$

However, these default values are valid only for stress states where the Eq. (15) is satisfied.

$$|\sigma_h| \leq \sqrt{3} f_c \tag{15}$$

$$\sigma_h = \frac{1}{3}(\sigma_{xp} + \sigma_{yp} + \sigma_{zp}) \tag{16}$$

Both the function F and failure surface S are expressed in terms of principal stresses denoted as $\sigma_1, \sigma_2$ and $\sigma_3$ as:

$$\sigma_1 = \max(\sigma_{xp}, \sigma_{yp}, \sigma_{zp}) \tag{17}$$

$$\sigma_1 = \min(\sigma_{xp}, \sigma_{yp}, \sigma_{zp}) \tag{18}$$

$$\sigma_1 \geq \sigma_2 \geq \sigma_3 \tag{19}$$

The failure of concrete is categorized into four domains:

1. $0 \geq \sigma_1 \geq \sigma_2 \geq \sigma_3$ (compression-compression-compression)

2. $\sigma_1 \geq 0 \geq \sigma_2 \geq \sigma_3$ (tensile-compression-compression)
3. $\sigma_1 \geq \sigma_2 \geq 0 \geq \sigma_3$ (tensile-tensile-compression)
4. $\sigma_1 \geq \sigma_2 \geq \sigma_3 \geq 0$ (tensile-tensile-tensile)

In each domain, independent functions describe $F$ and the failure surface $S$. These functions are described in detail below for each domain [29].

### 4.1. Compression-compression-compression domain

In the compression-compression-compression regime, $F$ and $S$ are defined as:

$$F = \frac{1}{\sqrt{15}}\left[(\sigma_1-\sigma_2)^2 + (\sigma_2-\sigma_3)^2 + (\sigma_3-\sigma_1)^2\right]^{\frac{1}{2}} \tag{20}$$

$$S = \frac{2r_2(r_2^2-r_1^2)\cos\eta + r_2(2r_1-r_2)\left[4(r_2^2-r_1^2)\cos^2\eta + 5r_1^2 - 4r_1 r_2\right]^{\frac{1}{2}}}{4(r_2^2-r_1^2)\cos^2\eta + (r_{2-}2r_1)^2} \tag{21}$$

in which:

$$\cos\eta = \frac{2\sigma_1-\sigma_2-\sigma_3}{\sqrt{2}\left[(\sigma_1-\sigma_2)^2 + (\sigma_2-\sigma_3)^2 + (\sigma_3-\sigma_1)^2\right]^{\frac{1}{2}}} \tag{22}$$

$$r_1 = a_0 + a_1\xi + a_2\xi^2 \tag{23}$$

$$r_2 = b_0 + b_1\xi + b_2\xi^2 \tag{24}$$

$$\xi = \frac{\sigma_h}{f_c} \tag{25}$$

$\sigma_h$ is defined by Eq. (16) and coefficients $a_0, a_1, a_2, b_0, b_1$ and $b_2$ are determined through solution of the simultaneous equations:

$$\left\{\begin{array}{l}\dfrac{F}{f_c}(\sigma_1=f_t,\sigma_2=\sigma_3=0)\\[4pt]\dfrac{F}{f_c}(\sigma_1=0,\sigma_2=\sigma_3=-f_{cb})\\[4pt]\dfrac{F}{f_c}(\sigma_1=-\sigma_h^a,\sigma_2=\sigma_3=-\sigma_h^a-f_1)\end{array}\right\}=\begin{bmatrix}1 & \xi_t & \xi_t^2\\ 1 & \xi_{cb} & \xi_{cb}^2\\ 1 & \xi_1 & \xi_1^2\end{bmatrix}\begin{Bmatrix}a_0\\ a_1\\ a_2\end{Bmatrix} \qquad (26)$$

$$\left\{\begin{array}{l}\dfrac{F}{f_c}(\sigma_1=\sigma_2=0,\ \sigma_3=-f_c)\\[4pt]\dfrac{F}{f_c}(\sigma_1=\sigma_2=-\sigma_h^a,\ \sigma_3=-\sigma_h^a-f_2)\\[4pt]0\end{array}\right\}=\begin{bmatrix}1 & -\dfrac{1}{3} & \dfrac{1}{9}\\ 1 & \xi_2 & \xi_2^2\\ 1 & \xi_0 & \xi_0^2\end{bmatrix}\begin{Bmatrix}b_0\\ b_1\\ b_2\end{Bmatrix} \qquad (27)$$

in which,

$$\xi_t=\dfrac{f_t}{3f_c},\quad \xi_{cb}=\dfrac{2f_{cb}}{3f_c},\quad \xi_1=-\dfrac{\sigma_h^a}{f_c}-\dfrac{2f_1}{3f_c},\quad \xi_2=-\dfrac{\sigma_h^a}{f_c}-\dfrac{f_2}{3f_c} \qquad (28)$$

and $\xi_0$ is the positive root of the equation:

$$r_2(\xi_0)=a_0+a_1\xi_0+a_2\xi_0^2=0 \qquad (29)$$

Since the failure surface must be remain convex, the ratio $r_1/r_2$ is restricted to the range

$$0.5<r_1/r_2<1.25 \qquad (30)$$

Although the upper bound is not considered to be restrictive since $r_1/r_2<1$ for most materials [29]. Also the coefficients $a_0,a_1,a_2,b_0,b_1$ and $b_2$ must satisfy the following conditions:

$$a_0>0, a_1\leq 0, a_2\leq 0 \qquad (31)$$
$$b_0>0, b_1\leq 0, b_2\leq 0 \qquad (32)$$

### 4.2. Tension-comparison-comparison domain

In the tension-comparison-comparison regime, $F$ and $S$ are defined as:

$$F = \frac{1}{\sqrt{15}}\left[(\sigma_2 - \sigma_3)^2 + \sigma_2^2 + \sigma_3^2\right]^{\frac{1}{2}} \quad (33)$$

$$S = \left(1 - \frac{\sigma_1}{f_t}\right)\frac{2p_2(p_2^2 - p_1^2)\cos\eta + p_2(2p_1 - p_2)\left[4(p_2^2 - p_1^2)\cos^2\eta + 5p_1^2 - 4p_1p_2\right]^{\frac{1}{2}}}{4(p_2^2 - p_1^2)\cos^2\eta + (p_2 - 2p_1)^2} \quad (34)$$

where $\cos\eta$ is defined by Eq. (22) and we have:

$$p_1 = a_0 + a_1\chi + a_2\chi^2 \quad (35)$$

$$p_2 = b_0 + b_1\chi + b_2\chi^2 \quad (36)$$

The coefficients $a_0, a_1, a_2, b_0, b_1$ and $b_2$ are defined by Eqs. (26) and (27) while

$$\chi = \frac{1}{3}(\sigma_2 + \sigma_3) \quad (37)$$

### 4.3. Tension-tension-compression domain

For this domain, F and S are defined as:

$$F = \sigma_i; \quad i = 1, 2 \quad (38)$$

$$S = \frac{f_t}{f_c}\left(1 + \frac{\sigma_3}{f_c}\right); \quad i = 1, 2 \quad (39)$$

### 4.4. Tension-tension-tension domain

For the last domain, we have:

$$F = \sigma_i; \quad i = 1, 2, 3 \quad (40)$$

$$S = \frac{f_t}{f_c}; \quad (41)$$

The other conditions are as defined for the first domain.

### 4.5. Constraints

In this optimization problem, various constraints including geometrical and stability constraints are considered. To ensure that upstream and downstream faces of the dam do not pass through each other Eq. (42) should be satisfied. Besides, for constructing facilities and having smooth cantilevers over the height of the dam, the slope of overhang at the upstream and downstream faces of the dam should satisfy Eq. (43):

$$r_{di} \leq r_{ui} \Rightarrow \frac{r_{di}}{r_{ui}} - 1 \leq 0, \; i = 1, 2, ..6 \tag{42}$$

$$\gamma \leq \gamma_{alw} \Rightarrow \frac{\gamma}{\gamma_{alw}} - 1 \leq 0 \tag{43}$$

where $r_{di}$ and $r_{ui}$ are the radius of curvature at *ith* level in the z direction; $\gamma$ is the slope of overhang at the downstream and upstream faces of dam. $\gamma_{alw}$ is allowable absolute value of aforementioned parameter. To ensure the sliding stability of the dam, following equation should be satisfied

$$\varphi_l \leq \varphi \leq \varphi_u \tag{44}$$

where $\varphi$ is the central angle of arch dam at *ith* level in z direction. $\varphi_l$ and $\varphi_u$ are the allowable lower and upper bounds of the central angle. According to theory of thin-walled structures the optimum value for central angle that minimizes the dam volume, is $133°, 34'$ but it has many restrictions due to construction reasons therefore $\varphi$ can vary from $90°$ to $130°$ throughout the dam height considering practical purposes.

## 5. Brief Review on The Single-objective Charged System Search

The Charged System Search (CSS) is a population-based search approach, which is inspired by the Coulomb law from electrostatics and the laws of motion from Newtonian mechanics [13]. In the CSS each agent (CP) is considered as a charged sphere with radius a, having a uniform volume charge density which can produce an electric force on the other CPs. The force magnitude for a CP located in the inside of the sphere is proportional to the separation distance between the CPs, while for a CP located outside the sphere it is inversely proportional to the square of the separation distance between the particles. The resultant forces or acceleration and the motion laws determine the new location of the CPs. In this stage, each CP moves in the direction of the resultant forces and its previous velocity. From optimization point of view, this process provides a good balancing between the exploration and the exploitation paradigms of the algorithm, which can considerably improve the efficiency of the algorithm [13]. Three essential concepts including self-adaptation step, cooperation step, and competition step, are considered in this algorithm. Moving towards good CPs provides the self-adaptation step. Cooperating CPs to determine the resultant force acting on each CP supplies the cooperation step and having larger force for a good CP, comparing a bad one, and saving good CPs in the Charged Memory (CM) provide the competition step. Application of the CSS method for solving many benchmark and engineering problems shows that it outperforms evolutionary algorithms, and comparison of the results demonstrates the efficiency of the present algorithm [13].

## 6. Multi-objective Charged System Search

Some modifications and additional steps are considered in the single-objective charged system search (CSS) algorithm for presenting a multi-objective CSS method (MoCSS). The pseudo-code for the MoCSS algorithm can be summarized as follows:

**Step 1:** *Initialization*. Similar to the standard CSS, the initial positions of CPs are determined randomly in the search space and the initial velocities of CPs are assumed to be zero. The values of the fitness functions for the CPs are determined and the CPs are ranked based on Pareto dominance criteria:

$$\forall i \neq j \qquad \exists fit_i \mid fit_i(\mathbf{X}A) \leq fit_i(\mathbf{X}B) \wedge fit_j(\mathbf{X}A) \leq fit_j(\mathbf{X}B) \tag{45}$$

**Step 2:** *definition of CM*. All CPs belong to first rank (according to the Pareto dominance criteria) are stored in a memory, so called charged memory (CM).

**Step 3:** *Determination of forces on CPs*. The force vector is calculated for each CP as

$$\mathbf{F}_j = \sum_{i,i\neq j}\left(\frac{q_i}{a^3}r_{ij}\cdot i_1 + \frac{q_i}{r_{ij}^2}\cdot i_2\right)ar_{ij}\,p_{ij}(\mathbf{X}_i - \mathbf{X}_j)\begin{cases} j=1,2,\ldots,N \\ i_1=1, i_2=0 \Leftrightarrow r_{ij} < a \\ i_1=0, i_2=1 \Leftrightarrow r_{ij} \geq a \end{cases} \quad (46)$$

where $\mathbf{F}_j$ is the resultant force acting on the $j$th CP; $N$ is the number of CPs. The magnitude of charge for each CP ($q_i$) is defined considering the quality of its solution as

$$q_i = \prod_{k=1}^{N}\frac{fit_k(i) - fitworst_k}{fitbest_k - fitworst_k}, \quad i = 1,2,\ldots,N \quad (47)$$

where operator $\Pi$ is used for repeated multiplications. $fitbest_k$ and $fitworst_k$ are the best and the worst fitness of the $k$th objective function for all CPs, respectively; $fit_k(i)$ represents the fitness of the agent $i$. $N$ is the number of objective functions. $P_{ij}$ is the probability of moving each CP towards the others and is obtained using the following function:

$$p_{ij} = \begin{cases} 1 & rank(i) > rank(j) \\ 0 & otherwise \end{cases} \quad (48)$$

In Eq. (46), $ar_{ij}$ indicates the kind of force and is defined as

$$ar_{ij} = \begin{cases} +1 & rand < 0.8 \\ -1 & otherwise \end{cases} \quad (49)$$

where $rand$ represents a random number.

**Step 4:** *Solution construction.* Each CP moves to the new position and the new velocity is calculated as:

$$\mathbf{X}_{j,new} = rand_{j1}\cdot k_a \cdot \mathbf{F}_j + rand_{j2}\cdot k_v \cdot \mathbf{V}_{j,old} + \mathbf{X}_{j,old} \quad (50)$$

$$\mathbf{V}_{j,new} = \mathbf{X}_{j,new} - \mathbf{X}_{j,old} \quad (51)$$

where $k_a$ is the acceleration coefficient; $k_v$ is the velocity coefficient to control the influence of the previous velocity; and $rand_{j1}$ and $rand_{j2}$ are two random numbers uniformly distributed in the range (0,1).

**Step 5:** *Position updating process.* If a new CP exits from the allowable search space, a harmony search-based handling approach is used to correct its position. According to this method, any variable of each solution that violates its corresponding boundary can be regenerated from CM as Fig. 5. In this figure, "w.p." is the abbreviation of "with the probability"; CMCR (the Charge Memory Considering Rate) varying between 0 and 1 sets the rate of selecting a value in the new vector from historic values stored in CM, and (1–CMCR) sets the rate of randomly choosing one value from a possible range of values. The value (1–PAR) sets the rate of doing nothing, and PAR sets the rate of choosing a value from neighboring the best CP.

**Step 6:** *CM updating.* All new CPs and available ones in the CM are ranked using Eq. (45) and then the first Pareto are saved in the CM then the agents with small objective distances will be deleted from the memory. This distance is defined as [30, 31]:

$$d_{ij} = \sqrt{\sum_{k=1}^{N}\left(u_k\left(fit_k\left(X_i\right)-fit_k\left(X_j\right)\right)\right)^2} \tag{52}$$

in which, $d_{ij}$ is the objective distance between solution $X_i$ and $X_j$ and $u_k$ is a constant chosen to make all $u_k fit_k$ close to each other. To determine the appropriate value of $u_k$, first we set all $u_k = 1$. Then by optimizing the problem, some near optimal solutions and $fitworst_k$ for each objective are obtained. Eventually, $u_k$ is obtained in the following way:

$$u_1 = \alpha > 0, u_k = u_{k-1} * \frac{fitworst_k}{fitworst_{k-1}} \quad for\, fit_k, \quad k > 1 \tag{53}$$

**Step 7:** *Termination criterion control.* Steps 3-6 are repeated until a termination criterion is satisfied.

## 7. Multi Criteria Decision Making

The output of a multi objective optimization problem is finding a Pareto set which meets the decision maker's (DM) preferences. Since the Pareto optimal solutions cannot be ranked globally, the decision maker needs to define some extra information to reach his desired answer. To date, various approaches have been invented for decision-making [32]. Here, we use the multi-criteria tournament decision-making method (MTDM) [33] because it is a simple method, which uses the preferences of decision maker in the form of criterion importance weight to rank solutions from best to worst. In this method, the DM global interests can be defined as a function R. To fulfill this aim, tournament function $T_i(a,A)$, should be defined as following:

$$T_i(a,A) = \sum_{\forall b \in A, a \neq b} \frac{t_i(a,b)}{(|A|-1)}, \tag{54}$$

where:

$$t_i(a,b) = \begin{cases} 1 & \text{if } fit_i(b) - fit_i(a) > 0 \\ 0 & \text{otherwise} \end{cases} \tag{55}$$

This function counts the ratio of times alternative $a$ wins the tournament against each other $b$ solution from $A$ in which, $a$ is a non-dominated point in the objective space.

The function $R$ is then defined in order to gather all individual criterion and their appropriate priority weights, $w_i$, into the global ranking function:

$$R(a) = \left( \prod_{i=1}^{N} T_i(a,A)^{w_i} \right)^{\frac{1}{N}} \tag{56}$$

in which $N$ is the number of objective functions and $w_i$ should be specified by DM under following conditions:

$$w_i > 0 \quad \text{and} \quad \sum_{i=1}^{n} w_i = 1. \tag{57}$$

## 8. Numerical Investigation

In order to show the performance of proposed methodology, the Morrow point arch dam as a well-known and real-world structure is selected as the case study in this paper. The MoCSS are coded in MATLAB software, modeling and analyzing of the arch dam are handled using a combination of parallel working MATLAB and Ansys Parametric Design Language (MATLAB-APDL) codes. Reservoir is supposed to be full and the dam-reservoir interaction subjected to seismic loading is taken to account in this example. In order to construct the dam geometry, six controlling levels are considered so the dam can be modeled using twenty shape design variables:

$$\mathbf{X} = \{\gamma, \beta, t_{c1}, t_{c2}, t_{c3}, t_{c4}, t_{c5}, t_{c6}, r_{u1}, r_{u2}, r_{u3}, r_{u4}, r_{u5}, r_{u6}, r_{d1}, r_{d2}, r_{d3}, r_{d4}, r_{d5}, r_{d6}\} \tag{58}$$

In optimization process, we need the lower and upper bounds of design variables; these can be obtained using classic design approaches, as [34]:

$$
\begin{aligned}
&0 \leq \gamma \leq 0.3 \qquad\qquad 0.5 \leq \beta \leq 1 \\
&3m \leq t_{c1} \leq 10m \qquad 104m \leq r_{u1} \leq 135m \qquad 104m \leq r_{d1} \leq 135m \\
&5m \leq t_{c2} \leq 14m \qquad 91m \leq r_{u2} \leq 118m \qquad 91m \leq r_{d2} \leq 118m \\
&7m \leq t_{c3} \leq 19m \qquad 78m \leq r_{u3} \leq 101m \qquad 78m \leq r_{d3} \leq 101m \\
&9m \leq t_{c4} \leq 23m \qquad 65m \leq r_{u4} \leq 85m \qquad 65m \leq r_{d4} \leq 85m \\
&11m \leq t_{c5} \leq 26m \qquad 52m \leq r_{u5} \leq 68m \qquad 52m \leq r_{d5} \leq 68m \\
&12m \leq t_{c6} \leq 31m \qquad 39m \leq r_{u5} \leq 51m \qquad 39m \leq r_{d6} \leq 51m
\end{aligned} \tag{59}
$$

Since selection of any design earthquake will not affect the optimization process of the proposed methodology, in this research the N–S record of 1940 El Centro earthquake is selected to apply to the arch dam-reservoir system in the upstream–downstream direction [35] as shown in Fig.6.

In addition to the proposed MoCSS, the multi objective particle swarm optimization (MoPSO) [36] and a fast elitist non-dominated sorting genetic algorithm (NSGA-II) [30] are utilized as the well-known multi-objective optimization methods. The parameters of these methods are presented in Table 3 in which $N_{var}$ is the number of decision variables. In Table 4, the extreme values obtained by the MoCSS are compared to those of the MoPSO and NSGA-II. It should be noted that the maximum allowable values for the dam's volume and

Willam-Warnke failure criterion is limited to $3.4 \times 10^5$ and 1.3, respectively. According to this table, for those points which $fit_2$ is more important than $fit_1$, the MoCSS yields 11.5% and 4.1% better than the MoPSO and NSGA-II, respectively. Since the accepted results should be negative for $fit_2$, so Table 4 presents best extreme value with negative $fit_2$, as well. The MoCSS can find an acceptable design with $2.246 \times 10^5$ m³ while the best acceptable results for the MoPSO and NSGA-II are $2.522 \times 10^5$ m³ and $2.748 \times 10^5$ m³, respectively. Clearly, the best result is obtained by the MoCSS. The values for this design are presented in Table 5. The Pareto fronts of three methods are presented in Fig. 7. According to this figure, the results of the MoCSS are obviously better the two other methods. Fig. 8 shows the optimum shape of dam with the smallest volume obtained by the MoCSS. Since stress state of the dam body varies in every time step of the earthquake for each node, to obtain a single number for comparison we chose the maximum value of failure criterion's time history. It is obvious that if maximum value is negative (the negative values of failure criteria are desirable), rest of values will be satisfied. Fig. 9 shows the failure criterion time history for six levels described in Fig. 8. The maximum value is about -0.004 which occurred in level 4 at time 4.34 second of the earthquake and therefore the optimum design presented in Table 5 is acceptable.

For the process of decision making and finding the best solution, the DM should notify their preferences by considering all the information integrated in the Pareto front. Before this process, we omit the unacceptable results (with positive $fit_2$) and then the DM is performed. In this example, in order to show the wide range of possible solutions, five different scenarios are considered. Fig. 10 shows the selected solutions obtained via the MoCSS corresponding to each considered scenario. The numerical values of these scenarios for all algorithms are presented in Table 6. As it can be seen, the all good results are found by the MoCSS. It should be noted that these scenarios become important when some different designs become necessary. In the other words, some designs with different safety factors should be studied considering the dam position and other environment conditions.

## 9. Conclusions

In this study, shape optimization of arch dams including the dam-reservoir interaction is formulated as a multi-objective problem. The EL Centro N-S component of Imperial Valley earthquake in 1940 is chosen as a ground motion. The double curvature arch dam geometry is modeled using twenty shape parameters in six controlling levels across the dam height. A comprehensive numerical model, considering all boundary conditions, using APDL is developed and verified. The optimization problem is formulated as two conflicting objectives including the concrete volume and Willam-Warnke failure criterion as a function of principal stresses in the dam body. Besides, some geometry and stability constraints are involved.

To solve this new problem, a new methodology, so-called Multi-objective charged system search (MoCSS) is proposed. To fulfill this aim, we developed an enhanced parallel-working APDL-MATLAB code. In contrast with previous works, the presented approach yields a Pareto solution instead of a single design in which decision maker can choose the best solution to meet his preferences. A multi-criteria tournament decision-making method is employed herein. In order to examine the effectiveness of proposed methodology, Morrow Point arch dam optimization is performed. The resulted Pareto fronts of MoCSS as well as extreme values are obtained and compared with two other approaches: NSGA-II and MoPSO. Besides, selected solutions by MTDM for five scenarios are obtained for MoCSS and compared to those of two methods. According to the results, proposed method is a reliable tool and can be considered as a suitable alternative to conventional methods of designing arch dams.

**Conflict of interest**

On behalf of all authors, the corresponding author states that there is no conflict of interest.

List of figures:

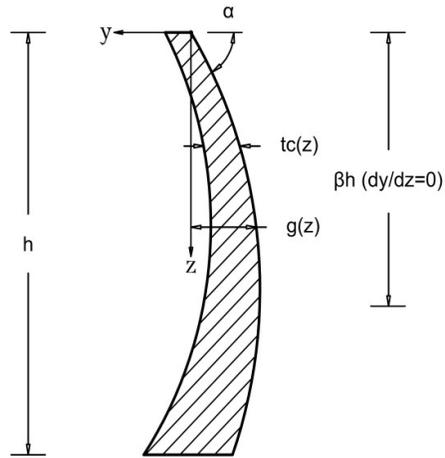

Fig. 1. Crown Cantilever profile of the arch dam

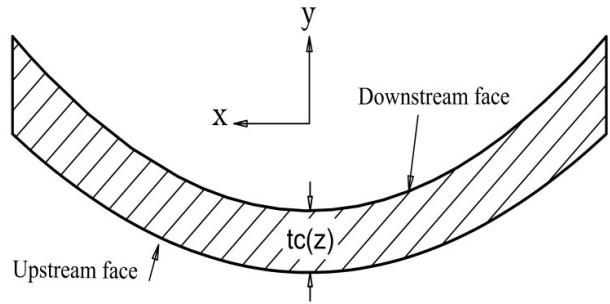

Fig. 2. Parabolic shape of an elevation of the arch dam

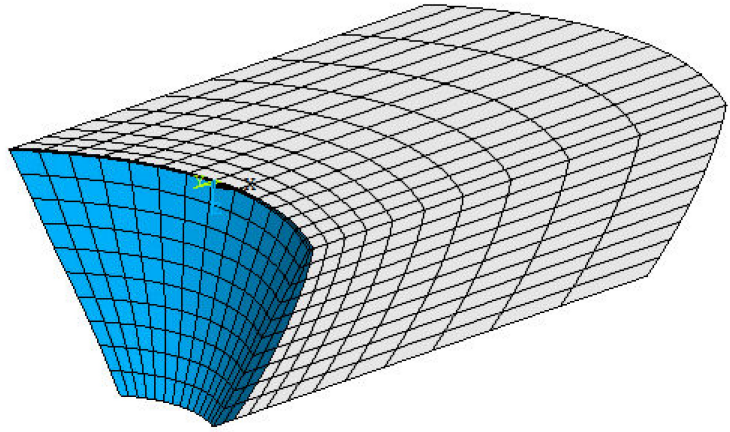

Fig. 3. Finite Element Model of the Morrow Point Dam

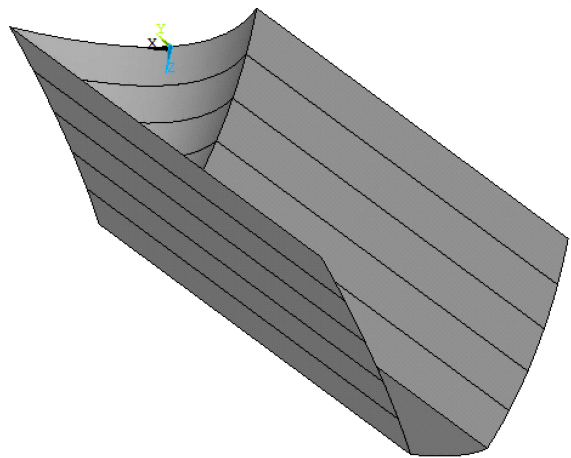

Fig. 4. Fluid-Structure Interaction Surfaces considered in the FEM

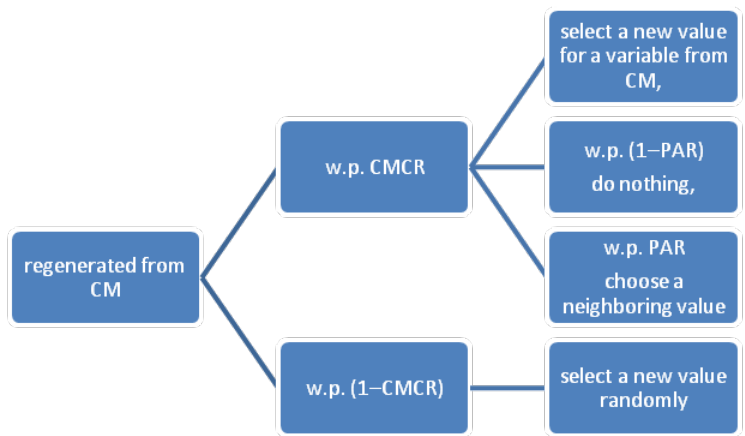

Fig. 5. Regenerating violated solutions from CM

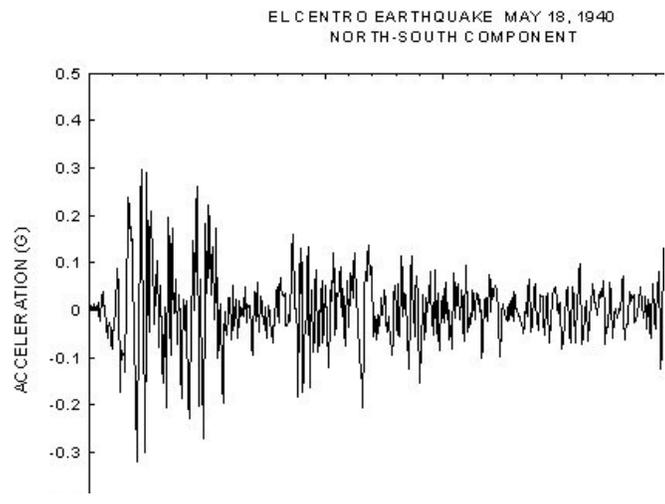

Fig. 6. N–S record of 1940 El Centro earthquake

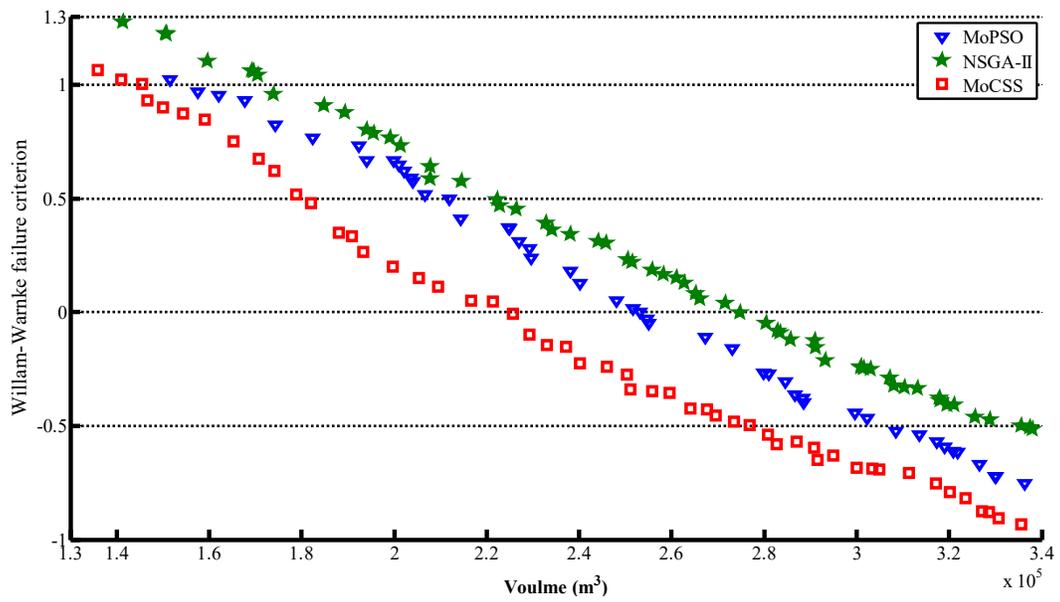

Fig. 7. Pareto fronts of three methods

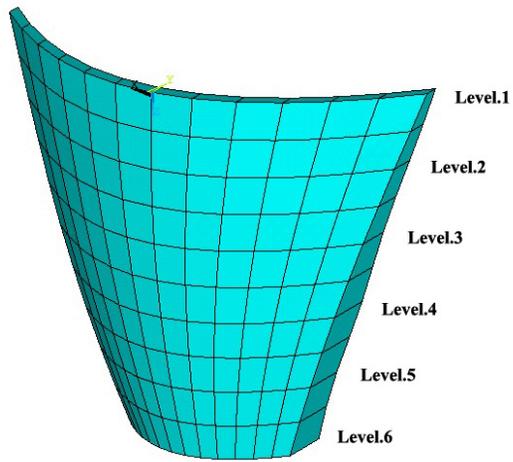

Fig. 8. Optimum shape of the Morrow Point dam resulted by the MoCSS

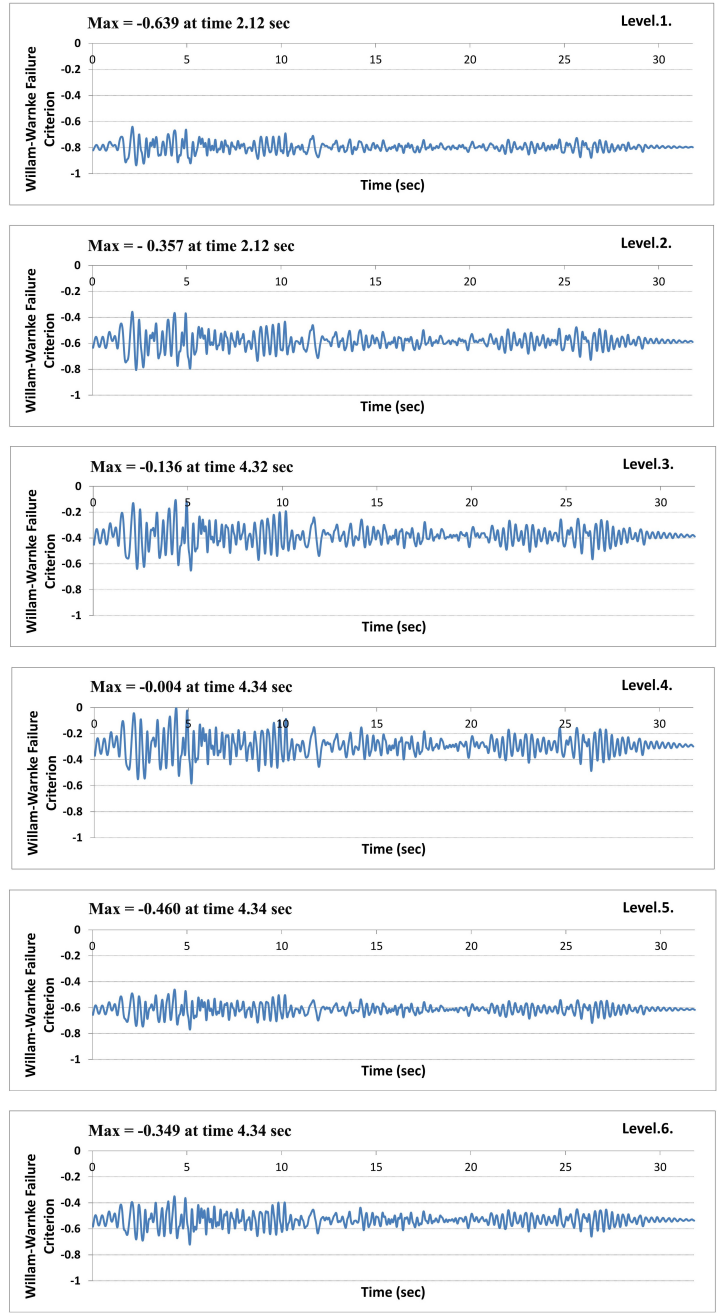

Fig.9. Failure criteria time history at different levels (Levels.1 to 6)

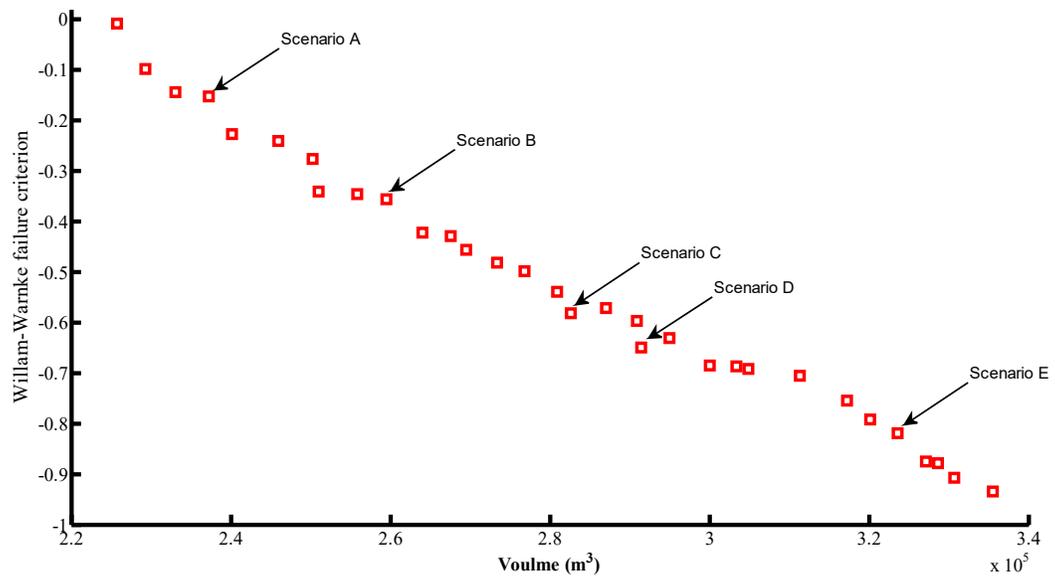

Fig. 10. Best solutions according to the five different scenarios

Table 1: Material Property of the Morrow Point Dam

| Material | Property | Value | unit |
|---|---|---|---|
| **Mass Concrete** | Modulus of Elasticity | 27.579 | GP |
| | Poisson's Ratio | 0.2 | - |
| | Mass Density | 2483 | Kg/m$^3$ |
| | Uniaxial compressive strength of concrete ($f_c$) | 30 | MPa |
| | Uniaxial tensile strength of concrete ($f_t$) | 1.5 | MPa |
| **Water** | Bulk Modulus | 2.15 | GP |
| | Mass Density | 1000 | Kg/m3 |
| | Velocity of Pressure Waves | 1438.66 | m/s |
| | Wave reflection coefficient | 1 | - |

Table 2: Comparison of the natural frequencies from the literature with FEM

| mode | Tan and Chopra [25] | | Duron and Hall [26] | | Present work | |
|---|---|---|---|---|---|---|
| | empty | Full | Full | | empty | full |
| | FEA | | FEA | experimental | FEA | |
| 1 | 4.27 | 2.82 | 3.05 | 2.95 | 4.2897 | 2.9607 |
| 2 | - | - | - | 3.3 | - | 3.3046 |
| 3 | - | - | 4.21 | 3.95 | - | 4.2425 |
| 4 | - | - | 5.96 | 5.4 | - | 4.9522 |
| 5 | - | - | - | 6.21 | - | 6.1156 |

Table 3: Parameters of the MoCSS method

| | | |
|---|---|---|
| **MoCSS** | Number of agents | 100 |
| | Maximum number of iterations | 200 |
| | $k_a$ | 2 |
| | $k_v$ | 2 |
| | $\alpha$ | 1 |
| | $k_t$ | 0.75 |
| **MoPSO** | Swarm size | 100 |
| | Archive size | 100 |
| | Mutation rate | 0.5 |
| | Number of divisions for the adaptive grid | 30 |
| **NSGA-II** | Population size | 100 |
| | Crossover probability ($p_c$) | 0.9 |
| | Mutation rate | $1/N_{var}$ |
| | Crossover distribution rate ($\eta_c$) | 20 |
| | Mutation distribution rate ($\eta_m$) | 20 |

Table 4: Comparison of the results for morrow point arch dam optimization problem

| Optimization method | MoCSS (present work) | NSGA-II | MoPSO |
|---|---|---|---|
| **Obtained extreme values** | [1.0650, 1.359e5] | [1.2740, 1.415e5] | [1.023, 1.516e5] |
| | [-0.9339, 3.355e5] | [-0.5714, 3.380e5] | [-0.7519, 3.364e5] |

Table 5: The best design obtained by the MoCSS

| Variable | Value | Variable | Value | Variable | Value |
| --- | --- | --- | --- | --- | --- |
| $\gamma(m/m)$ | 0.201 | $\beta(m/m)$ | 0.516 | $rd_1(m)$ | 109.716 |
| $tc_1(m)$ | 4.852 | $ru_1(m)$ | 110.637 | $rd_2(m)$ | 92.719 |
| $tc_2(m)$ | 8.974 | $ru_2(m)$ | 93.582 | $rd_3(m)$ | 79.562 |
| $tc_3(m)$ | 11.972 | $ru_3(m)$ | 80.408 | $rd_4(m)$ | 66.341 |
| $tc_4(m)$ | 16.298 | $ru_4(m)$ | 67.690 | $rd_5(m)$ | 54.418 |
| $tc_5(m)$ | 15.883 | $ru_5(m)$ | 55.084 | $rd_6(m)$ | 39.995 |
| $tc_6(m)$ | 16.891 | $ru_6(m)$ | 41.713 | | |

Table 6: Different possible scenarios for the morrow point arch dam with corresponding solutions

| Scenario | Importance of criteria | Possible priority weights | Selected solution by MTDM | | | | | | | | |
|---|---|---|---|---|---|---|---|---|---|---|---|
| | | | MoCSS | | | NSGA-II | | | MoPSO | | |
| | | | $fit_1$ | $fit_2$ | $R_i$ | $fit_1$ | $fit_2$ | $R_i$ | $fit_1$ | $fit_2$ | $R_i$ |
| A | $C1 \gg C2$ | [0.9,0.1] | 237222.6 | -0.151 | 0.8499 | 283505.7 | -0.090 | 0.8498 | 279822.1 | -0.266 | 0.8500 |
| B | $C1 > C2$ | [0.7,0.3] | 259548.8 | -0.356 | 0.7367 | 301487.4 | -0.248 | 0.7368 | 313505.0 | -0.538 | 0.7440 |
| C | $C1 \approx C2$ | [0.5,0.5] | 282659.1 | -0.581 | 0.7071 | 317859.0 | -0.377 | 0.7069 | 336399.4 | -0.751 | 0.7177 |
| D | $C1 < C2$ | [0.3,0.7] | 291472.7 | -0.648 | 0.7367 | 337451.3 | -0.508 | 0.7368 | 358781.1 | -0.941 | 0.7367 |
| E | $C1 \ll C2$ | [0.1,0.9] | 327110.7 | -0.875 | 0.8499 | 384455.6 | -0.566 | 0.8498 | 380008.4 | -0.999 | 0.8500 |